\documentclass[12pt]{article}
\usepackage{amssymb}

\def\appendix{\par}  

\def\res{\upharpoonright}                                                                               
\def\om{\omega}
\def\al{\alpha}
\def\be{\beta}
\def\ga{\gamma}
\def\de{\delta}
\def\sq{\subseteq}
\def\ss{\setminus}
\def\sep{{\rm Sep}}
\def\emp{\emptyset}
\def\bPi{{\bf \Pi}}
\def\bSig{{\bf \Sigma}}
\def\bDel{{\bf \Delta}}
\def\proof{\par\noindent Proof\par\noindent}
\def\qed{\par\noindent QED\par}
\def\rmiff{\mbox{ iff }}
\def\rmand{\mbox{ and }}
\def\nembed{\not\hookrightarrow}
\def\reals{{\mathbb R}}
\def\rationals{{\mathbb Q}}
\def\adm{{\mathbb A}}
\def\ll{{G{\cal N}}}
\def\rank{{\mbox{ rank}}}

\newtheorem{theorem}{Theorem}
\newtheorem{lemma}[theorem]{Lemma}
\newtheorem{question}[theorem]{Question}
\newtheorem{prop}[theorem]{Proposition}

\begin{document}

\begin{center}
Half of an inseparable pair
\end{center}

\begin{flushright}
Arnold W. Miller\footnote{
Thanks to Jindrich Zapletal who organized the SEALS 
meeting at the University of Florida, Gainesville in March 2004
during which part of these results were obtained.
\par LaTeX2e - texed on \today
} 
\end{flushright}

\begin{quote}
Abstract: 
A classical theorem of Luzin is that the separation principle holds for the
$\bPi^0_\al$ sets but fails for the $\bSig^0_\al$ sets.  We
show that for every $\bSig^0_\al$ set $A$ which is not $\bPi^0_\al$
there exists a $\bSig^0_\al$ set $B$ which is disjoint from $A$ but cannot
be separated from A by a $\bDel^0_\al$ set $C$.  Assuming 
$\bPi^1_1$-determancy it follows from a theorem of Steel that a similar result
holds for  $\bPi^1_1$ sets.  On the other hand assuming V=L there is a
proper $\bPi^1_1$ set which is not half of a Borel inseparable pair. 
These results answer questions raised by F.Dashiell. 
\end{quote}

The separation principle is a classical property of point classes in
descriptive set theory.  For every countable ordinal $\al$ and every pair of 
disjoint sets $A,B\sq 2^\om$ in the multiplicative class  $\al$
($\bPi^0_\al$) there exists a set $C$ in ambiguous class $\al$ ($\bDel^0_\al$)
which separates them, i.e.,  $A\sq C$ and $C\cap B=\emp$.  It is also classical
result of Luzin 
that the separation principle must fail for the dual classes $\bSig^0_\al$.  
For proofs, see Kechris \cite{kechris} \S 22.  

For $\Gamma$ a class of subsets of $\om^\om$, define the
dual class $\widetilde{\Gamma}=\{\om^\om\ss A:A\in\Gamma\}$, 
$\Delta=\Gamma\cap \widetilde{\Gamma}$, and 

$\sep(\Gamma)\equiv \forall A,B\in \Gamma\;\; A\cap B=\emp\to
\exists C\in \Delta\;\; A\sq C \rmand A\cap B=\emp$.

\noindent $\Gamma$ is continuously closed iff for all continuous
$f:\om^\om\to \om^\om$ if $A\in \Gamma$ then $f^{-1}(A)\in\Gamma$.
 $\Gamma$ is nonselfdual iff $\Gamma\neq\widetilde{\Gamma}$.

Van Wesep and Steel
\cite{vanwesep} \cite{vanwesep2} \cite{steelsep} proved that 
for continuously closed nonselfdual $\Gamma$ in
the Borel subsets of $\om^\om$ that either 
($\neg \sep(\Gamma)$ and $\sep(\widetilde{\Gamma}))$ or
($\neg\sep(\widetilde{\Gamma})$ and $\sep(\Gamma)$), i.e.,
separation holds on one side and fails on the other.  This result 
is true for all continuously closed nonselfdual classes,
if the Axiom of Determinacy holds.

In Dashiell \cite{dashiell}, Luzin's theorem on the 
failure of separation for ${\bSig}^0_\al$ is used to
prove that the Banach space, ${\mathcal B}_\al$, of Baire class
$\al$-functions is not isomorphic to the space ${\mathcal B}_{\om_1}$ 
of Baire functions.  

The following Theorem settles a question raised by F. Dashiell.  He already
knew the result for ${\bSig}^0_1$ and ${\bSig}^0_2$.  It was also 
apparently\footnote{Unfortunately, my French is not very good, but I think this 
may be the question on the top of page 73 of \cite{luzin}, ``Un autre
probl\`eme \ldots '' and the last paragraph on page 76.} 
asked by Luzin \cite{luzin} in 1930.

\begin{theorem}\label{one}
Suppose $X$ is a Polish space and $A\sq X$ is
${\bSig}^0_\al$ but not ${\bPi}_\al^0$.  Then
there exists $A^*\sq X$ which is ${\bSig}^0_\al$ such that
$A\cap A^*=\emp$ but there does not exist a ${\bDel}_\al^0$ set $C$  
which separates $A$ and $A^*$, i.e., 
$A\sq C$ and $C\cap A^*=\emp$.
\end{theorem} 

\proof
For $\al=1$, if $A$ is any open set which is not closed, then
it cannot be separated from the interior of $X\ss A$.  So we
may assume $\al\geq 2$.  By Theorem 4 of Kunen-Miller \cite{km},
there exists a set $P\sq X$ such that $P$ is homeomorphic to a
closed subset of 
$2^\om$ and $P\cap A$ is ${\bSig}^0_\al\ss {\bDel}_\al^0$.
So without loss of generality we may assume $A\sq 2^\om$.

For subsets $B,C\sq 2^\om$ define $B\leq_W C$ (Wadge reducible) iff
there exists a continuous map $f:2^\om\to 2^\om$ such that
$f^{-1}(C)=B$.  Associated with Wadge reducibility is the Wadge game
whose payoff set is of roughly 
the same complexity as $B$ and $C$.  It follows
from Borel determinacy, see Martin \cite{martin},
that for every pair of Borel sets $B$ and $C$ that either $B\leq_W C$ or
$C\leq_W (2^\om\ss B)$, see for example Van Wesep \cite{vanwesep}.
It follows from this that for any $B\sq 2^\om$ which is ${\bSig}^0_\al$
we have that $B\leq_W A$, since otherwise $A\leq_W(2^\om\ss B)$ would
make $A$ a ${\bPi}^0_\al$ and hence ${\bDel}^0_\al$, which is
contrary to our assumption.

Now assume $\al=2$. Let $D ,D^*\sq 2^\om$ be countable dense and
disjoint.  Note that they are ${\bSig}^0_2$ sets which cannot be
separated, since dense ${\bPi}^0_2$, i.e., $G_\de$, sets must
intersect by the Baire Category Theorem.
Since $D \leq_W A$ there exists a continuous map
$f:2^\om\to 2^\om$ with $f^{-1}(A)=D $. 
Let $A^*=f(D^*)$.  Since it
is countable, $A^*$ is a ${\bSig}^0_2$ set.  It cannot be separated from
$A$, because if $C$ is a ${\bDel}^0_2$ with $A\sq C$ and
$A^*\cap C=\emp$, then $D \sq f^{-1}(C)$ and 
$D^*\sq f^{-1}(2^\om\ss C)$ would separate $D $ and $D^*$.

Now assume $\al>2$.
By a result of Harrington, see Steel \cite{steel} or Van Engelen, Miller,
Steel \cite{vms}, for any $B$ which is ${\bSig}^0_\al$ there exists
a one-to-one continuous map $f:2^\om\to 2^\om$ such that 
$f^{-1}(A)=B$.  By a classical theorem of descriptive set theory
(see Kechris \cite{kechris}) there exists disjoint $B,B^*\sq 2^\om$ 
${\bSig}^0_\al$ sets which cannot be separated by a 
${\bDel}^0_\al$ set.  Let $f$ be one-to-one and continuous with $f^{-1}(A)=B$.
Let $A^*=f(B^*)$.  Since $f$ is one-to-one, it is a homeomorphism onto
its range and hence $A^*$ is a ${\bSig}^0_\al$ set disjoint from
$A$.  The set $A^*$ cannot be separated from $A$ because 
the preimage of a separating
set would separate $B$ and $B^*$.

\qed

Dashiell's proof of Theorem \ref{one} for $\al=2$ is as follows. Suppose $X$
is a Polish space and $A\sq X$ is some $F_\sigma$ set which is not a
$G_\de$.   By Baire's theorem on functions of the first class, there exists
a closed  $F \sq X$ on which the characteristic function of $A$ has no point of
continuity relative to $F$.  That is, both $A\cap F$ and $A\ss F$ are
dense in F.  Let $A^*$ be a countable dense set in $A\ss F$ (hence an
$F_\sigma$).  Clearly now $A$ and $A^*$ can not be separated by disjoint
$G_\de$ sets of X, because intersecting with $F$ would  give two dense
$G_\de$ subsets of the complete metric space $F$, which must meet. 

Dashiell pointed out that for a fixed countable ordinal $\al$ if we let $X_\al$
be the Stone space of the Boolean algebra of $\bDel_\al^0$ subsets of the
reals, then the cozero sets in $X_\al$ whose closures are not open (i.e., not
clopen) correspond to the proper $\bSig^0_\al$ sets.   Hence, by Theorem
\ref{one}, we know that every cozero set $A$ whose closure is not open has
an inseparable disjoint sibling, i.e., a cozero set $B$ disjoint from $A$ but
the closures of $A$ and  $B$ must meet.

Dashiell tells us that the question from \cite{dashiell} of whether 
${\cal B}_\al$ and ${\cal B}_\be$ can be isomorphic Banach spaces for some
$1<\al<\be<\om_1$  is still open.

Dashiell also raised the same question for the coanalytic sets, $\bPi^1_1$.
The classic result (see Kechris \cite{kechris} \S 34,35)
is that any pair of disjoint analytic sets ($\bSig_1^1$)
can be separated by a Borel set ($\bDel^1_1$), but separation
fails for $\bPi^1_1$.  Luzin proved this by applying the reduction
principle to a pair of doubly universal sets.

\begin{theorem}\label{two}
Suppose $\bPi^1_1$-determinacy holds, then for any $\bPi^1_1$ set $A$
in a Polish space $X$, if
$A$ is not $\bSig^1_1$, then there exists $A^*\sq X$ a $\bPi^1_1$ set 
disjoint from $A$
which cannot be separated from $A$ by a Borel set ($\bDel^1_1$).
\end{theorem}

\begin{theorem}\label{three}
Suppose $V=L$,  then there exists a $A\sq 2^\om$ and $\bPi^1_1$ set
which is not $\bSig^1_1$ with the
property that for any $B\sq 2^\om$ a $\bPi^1_1$ set disjoint from $A$
there exists a Borel set $C$ with $A\sq C$ and $C\cap B=\emp$.
\end{theorem}

\proof
For Theorem \ref{two} note that
since there is a Borel bijection between $X$ and $2^\om$ we may assume
that $X=2^\om$.
Theorem \ref{two} is an immediate corollary of a Theorem of Steel
\cite{steel}, who showed that $\bPi^1_1$-determinacy implies that for
any two properly $\bPi_1^1$ subsets $A_1,A_2$ of $2^\om$ there
exists a Borel automorphism $f:2^\om\to 2^\om$ such that
$f(A_1)=A_2$.  Hence if we take $C,C^*\sq 2^\om$ to be any disjoint
pair of $\bPi^1_1$ sets which are not Borel separable and
$f:2^\om\to 2^\om$ a Borel automorphism with $f(A)=C$, then
$f^{-1}(C^*)=A^*$ will be the required set.  

For Theorem \ref{three} we use for $A$ the self-constructible reals
studied by Guaspari, Kechris, and Sacks, see Kechris \cite{kechris2}
\S 2, 
where the self-constructible reals $A$ are denoted
${\mathcal C}_1$.

Define 
$$A=\{x\in 2^\om\; :\; x\in L_{\om_1^x}\}$$
where $\om_1^x$ is the least ordinal which is not the order type
of a relation recursive in $x$.  It is also the least ordinal $\al$
such that $L_\al[x]$ is an admissible set.  Suppose that $B$ is
a $\bPi^1_1$ set disjoint from $A$.  Then we may assume
that $B$ is $\Pi^1_1(x)$ for some $x\in A$ since by Kechris \cite{kechris2}
2A, every real in $L$ is recursive in some $x\in A$. 

The following Lemma is a relativized version
of Sacks \cite{sacks} III Lemma 9.3 p. 82.  For the convenience
of the reader we give a proof.  We are cheating a little bit,
by using the Addison-Kondo Theorem since Sack's uses his lemma
to deduce this result.  But of course it is OK since Addison-Kondo
has other proofs.   

Let $\ga<\om_1^x$ be the least ordinal so that $x\in L_\ga$.  For any $y\in
2^\om$ define $\ga^+(y)$ to be the least $\al>\ga$ such that $L_\al[y]$ is an
admissible set. 

\begin{lemma}\label{mainlemma}
For any $C\sq 2^\om$ a nonempty 
$\Pi^1_1(x)$ set there exists $y\in C$ such that
$y\in L_{\ga^+(y)}$.
\end{lemma}
\proof

Recall that a binary relation $(X,R)$ is well-founded iff
every nonempty subset of $X$ has an $R$-minimal element.  A
map $f:X\to $ Ordinals is called a rank function iff
$$\forall s,t\in X\;\; sRt \to f(s)<f(t).$$ 
Then $(X,R)$ is well-founded iff
it has a rank function on it.  For $(X,R)$ well-founded
the canonical rank function on $X$ is defined inductively
by 
 $$f(s)=\sup\{f(t)+1:tRs\}.$$
The range of the canonical
rank function is called the rank of $(X,R)$.
Furthermore, if $(X,R)\in \adm$ is a well-founded relation in
an admissible set $\adm$, then its rank and its canonical rank function 
are in $\adm$. See Barwise \cite{barwise} V.3.1 p.159.  

\bigskip\noindent
{\bf Claim 1}. Suppose $T\sq \de_1^{<\om}$ is a subtree, $T\in L_{\de_2}$ where
$\de_2>\om$ is a limit ordinal. For
each $s\in T$ define $T_s=\{t\in T\;:\; s\sq t\}$.  For each ordinal
$\al<\de_2$ if $\rank(T_s)=\al$ then the canonical rank function,
on $T_s$, i.e., $t\mapsto\rank(T_t)$ is an element of $L_{\de_2+\al+1}$.

\proof
Note that $(T\times\al)\in L_{\de_2}$ since $\al$ is small.
Fix $\al$ and $s\in T$ with $\rank(T_s)=\al$.
For each $\de<\de_1$ if $s\de\in T$ and 
$\rank(T_{s\de})=\beta$, then the canonical rank
function on $T_{s\de}$ is in $L_{\de_2+\be+1}\sq L_{\de_2+\al}$
and is uniformly definable from $T_{s\de}$, hence
the canonical rank function on $T_s$ is in 
$L_{\de_2+\al+1}$.
\qed

\bigskip\noindent
{\bf Claim 2}. Suppose $T$, $\de_1$ and $\de_2$ satisfy the hypothesis
of Claim 1.  For any ordinal $\al$ define
$$T(\al)=\{s\in T\;:\; \rank(T_s)<\al\}.$$
Then $T(\al)\in L_{\de_2+\al+1}$.

\proof
This follows from the previous claim since the canonical rank functions
are elements of $L_{\de_2+\al}$.
\qed

By the Addison-Kondo Theorem we may assume that $C$ is a
$\Pi^1_1(x)$ singleton, i.e. $C=\{y_0\}$.  

Now by standard arguments there exists a tree  
 $T\sq \cup_{n<\om} (\om^{n}\times 2^{n})$
which is recursive in $x$ such that for every $y\in 2^\om$
we have that
 $$y=y_0 \rmiff T\langle y\rangle
 =^{def}\{s\;:\; (s,y\res |s|)\in T\}\sq \om^{<\om} 
 \mbox{ is well-founded. }$$
Now since the tree $(T\langle{y_0}\rangle,\supset)$ is well-founded  
and it is an element of the admissible set $L_{\ga^+(y)}[y]$,
its rank $\de_0$ is strictly less than $\ga^+(y)$ and its canonical
rank function $R:T\langle{y_0}\rangle\to \de_0$ is in $L_{\ga^+(y)}[y]$.

Now define a tree 
$$T^*\sq \cup_{n<\om} (\de_0^{n}\times 2^{n})$$
which basically consists of attempts at a rank function 
into $\de_0$ for $T\langle{y_0}\rangle$. 
More formally, suppose $\{t_i:i<\om\}$ is a reasonable recursive listing
of $\om^{<\om}$, e.g., it should have the properties 
that $|s_i|\leq i$ and if $s_i\subset s_j$ then $i<j$.  

Define
$(r,s)\in T^*\cap(\de^{n}\times 2^{n})$ iff 
for each $i,j<n$ 

if $(t_i,s\res|t_i|),(t_j,s\res|t_j|)\in T$ and $t_i\subset t_j$ then
$r(j)<r(i)$.

\noindent Let $R^*:\om\to\de_0$ be the corresponding map to $R$, i.e.,
$$R^*(i)=
\left\{ 
\begin{array}{ll}
R(t_i) & \mbox{if }t_i\in T\langle{y_0}\rangle\\
0      & \mbox{otherwise}
\end{array}\right.$$
Note that $T^*$ is an element of $L_{\ga^+(y_0)}$
and $(y_0,R^*)$ is an infinite branch thru it.
We claim that $(y_0,R^*)$ is the lexicographically least
infinite branch thru $T^*$.  To see this, note that
if $(y,S)$ is an infinite branch in $T^*$, then $y=y_0$, since
$S$ will be a rank function for $T\langle y\rangle$, hence $T\langle y\rangle$ 
is well-founded and so $y=y_0$.  On the other hand $R$ assigns to
any $s\in T\langle{y_0}\rangle$ the smallest possible ordinal for
any rank function, and so $R^*$ will 
be lexicographically less than $S$.

Let
$$LF=\{\sigma\in T^*: \sigma \mbox{ is lexicographically left of }
(y_0,R^*)\}.$$
Then $(LF,\supset)$ is a well-founded relation and it is an element of
the admissible set $L_{\ga^+(y_0)}[y_0]$.  Hence its rank $\de_1$
is strictly smaller than $\ga^+(y_0)$.  By identifying the tree $T^*$ 
with a tree on $(\de_0+\de_0)^{<\om}$, i.e., by mapping
$(i,\al)\in 2\times\de_0$ to $\de_0\cdot i+\al$ we may apply 
Claim 2.  Hence the tree
$T^*\ss T^*(\de_1)$
and its leftmost branch $(y_0,R^*)$ (which is $\Delta_1$ in it)
are elements of $L_{\ga^+(y_0)}$.

Hence $y_0\in L_{\ga^+(y_0)}$ as was to be shown.
\qed

The relation
 $$\{(u,v):u\in\Delta^1_1(v)\}$$
is $\Pi^1_1$.  Hence the set
 $$C=\{y\in B: x\in\Delta^1_1(y)\}$$
is $\Pi_1^1(x)$.  If it is nonempty, then there exists $y\in C$
with $y\in L_{\ga^+(y)}$.  But since $x\in\Delta^1_1(y)$ we
know that $\om_1^y\geq \om_1^x>\ga$ hence $y\in L_{\om_1^y}$ which
contradicts $A\cap B=\emp$.  It follows that
$$B\sq \{y: x\notin \Delta^1_1(y)\}\sq \{y: \om_1^y < \ga\}$$
The second inclusion is true since every element of $L_{\om_1^y}$ is
in $\Delta_1^1(y)$.  It is well known that for any countable
$\ga$ the set $D=\{y\in 2^\om: \om_1^y < \ga\}$ is Borel. 
For example, a $\bSig^1_1$ definition and $\bPi_1^1$ definition are
given by:
\begin{enumerate}
 \item $y\in D$ iff
 there exists $\al<\ga$ such that $\forall e\in \om $ if $\{e\}^y$ is
 characteristic function of 
 a well-ordering $(\om,\leq_e^y)$, then $\;\;$ order-type$(\om,\leq_e^y)<\al$.
 \item $y\in D$ iff there does not exist $e\in\om$ and
 $f:(\om,\leq_e^y)\to(\ga,<)$ an isomorphism where $\{e\}^y$ is
 the characteristic function of the relation $(\om,\leq_e^y)$.
\end{enumerate}
But note that $D\cap A\sq L_\ga$ is countable and $B\sq D$,
so $A$ and $B$ can be separated
by a Borel set.
\qed

\begin{question}    
If every non Borel $\bPi^1_1$ set is half of an inseparable pair,
then is $\bPi_1^1$-determinacy true?   
\end{question}

See Harrington \cite{harrington}
for some properties of coanalytic sets which imply
$\bPi_1^1$-determinacy.

Cliff Weil raised the question of whether we
can get a large number of examples in Theorem \ref{three}, e.g., 

\begin{question}
Assuming V=L,
does there exist continuum many coanalytic sets which are pairwise non Borel
isomorphic and each of which is not half of an inseparable pair?
\end{question}

In Cenzer and
Mauldin \cite{cenzer} it is shown that assuming V=L there are continuum many 
coanalytic sets no two of which are Borel isomorphic.

\begin{center}
Separation for subsets of $\om$.
\end{center}

We could also consider the failure of separation for (lightface) classes of
subsets of $\om$.  Addison \cite{addison} shows that separation holds for
the class of $\Pi_n^0$ and fails for the class $\Sigma_n^0$ subsets of $\om$.  
However, not every proper $\Sigma_1^0$ subset of $\om$ is half of an
inseparable pair.   A set $A\sq\om$ is simple iff it is 
recursively enumerable (equivalently $\Sigma^0_1$),
coinfinite, but its complement does not contain an infinite recursively
enumerable subset.  Simple sets were first constructed by Post \cite{post} (or
see Soare \cite{soare}), and clearly a simple set cannot be half of an 
inseparable pair.  We are not sure exactly which recursively enumerable sets
are half of inseparable pair, perhaps just the complete ones.

Post also showed that a subset of $\om$ is $\Sigma^0_{n+1}$ 
iff it is
$\Sigma^0_1(0^{(n)})$ (see Soare\cite{soare} IV 2.2).  By relativizing his
construction of a simple set to the oracle $0^{(n)}$ we get a properly
$\Sigma^0_{n+1}$ subset of $\om$ which is not half of an inseparable pair.

Similarly, separation holds for the class of $\Sigma^1_1$ subsets of
$\om$ and fails for $\Pi^1_1$.  A proof analogous to the simple set 
type construction will give a proper $\Pi^1_1$ subset of $\om$
which is not half of an inseparable pair (see the proof of
Sacks \cite{sacks} VI Theorem 2.1 or 2.4).   
Another ``natural'' example of
such a $\Pi^1_1$-set can be given as follows.  Let
$(\om,\preceq)$ be a recursive linear ordering whose well-ordered initial
segment is isomorphic to $\om_1^{CK}$, the first non recursive ordinal.
The existence of such a linear ordering is due to Feferman \cite{feferman}
or perhaps Harrison \cite{harrison} see also Ash and Knight
\cite{ash} 8.11.  Now let
$A$ be the initial well-ordered segment of $\preceq$, i.e., 
 $$A=\{n\in\om: \{m: m\prec n\}\mbox{ is well-ordered by } \preceq\}.$$  
Then $A$ is a proper $\Pi^1_1$ set.  It cannot be half of an inseparable
pair because if $B\sq\om$ is $\Pi^1_1$ and disjoint from $A$ then there
must exists some $n_0\notin A$ such that $k\succeq n_0$ for every
$k\in B$.  Otherwise 
$$\om\ss A=\{m\in \om\::\; \exists k\in B\;\; k\preceq m\}$$
but $A$ is not a $\Delta^1_1$ set.

Another light-face question one might ask is the following.  Suppose
$A$ and $B$ are disjoint $\Pi^1_1$ subsets of $\om^\om$ which cannot
be separated by a $\Delta_1^1$-set, then can they be separated by
a $\bDel^1_1$-set?  Here is a counterexample.  Let $A,B\sq \om$ be
disjoint $\Pi^1_1$ sets which cannot be separated by
$\Delta^1_1$ subset of $\om$. 
Define $A^*=\{f\in\om^\om:f(0)\in A\}$ and
$B^*=\{f\in\om^\om:f(0)\in B\}$.  Then $A^*$ and $B^*$ are 
disjoint $\Pi^1_1$ which are clopen and hence separable 
by clopen sets.   But they cannot be separated
by a $\Delta^1_1$ subset of $\om^\om$.  Suppose
$C\sq\om^\om$ is $\Delta_1^1$ and $A^*\sq C$ and
$B^*\cap C=\emp$.  For each $n<\om$ let $x_n\in\om^\om$ be
the constant function $n$.  Then 
$$C^*=\{n<\om\;:\;x_n\in C\}$$
is a $\Delta_1^1$ set separating $A$ and $B$.

\begin{center}
Natural pairs of inseparable sets.
\end{center}

A number of authors have given natural examples of inseparable pairs
of $\bPi^1_1$ sets.

\medskip\noindent Luzin \cite{luzinbook} p.263 gives the following example.  
Let
$$\phi:\om^\om\times\om^\om\to\om^\om$$
be a Borel function such that for every $f:\om^\om\to\om^\om$
continuous there exists $x$ such that 
$\forall y\;\;\phi(x,y)=f(y)$.  Let
$$E=\{(x,z):\exists ! y\; \phi(x,y)=z\}$$
$$E_0=\{(x,z)\in E:\exists ! y\; \phi(x,y)=z \rmand y(0)\mbox{ is even }\}$$
$$E_1=\{(x,z)\in E:\exists ! y\; \phi(x,y)=z \rmand y(0)\mbox{ is odd }\}$$
Then $E_0$ and $E_1$ are disjoint inseparable $\bPi^1_1$ sets.

\medskip\noindent I wasn't able to decipher Novikov's example \cite{novikov}.

\medskip\noindent Sierpinski \cite{sier} gives the 
following pair of inseparable 
$\bPi^1_1$ sets.  Let $U\sq \reals^3$ be a universal $G_\de$ set
for subsets of the plane, i.e., $U$ is $G_\de$ and for every $G_\de$
set $V\sq \reals^2$ there exists an $x\in\reals$ with $U_x=V$.  Then

$S_1=\{(x,y)\;:\;\neg\exists z\; (x,y,z)\in U\}$

$S_2=\{(x,y)\;:\;\exists!\; z\; (x,y,z)\in U\}$

\noindent are a pair of inseparable $\bPi^1_1$ subsets of the plane.

\medskip\noindent Dellacherie and Meyer \cite{del} give the following pair
of inseparable $\bPi^1_1$ sets (or perhaps the analogous
families of trees):  Let
$LO$ be the space of linear orderings on $\om$ which
we can regard as a closed subspace of 
$P(\om\times\om)\equiv 2^{\om\times\om}$.
Let $WO\sq LO$ be the well-orderings.
For two linear orderings let $L_1\nembed L_2$ mean that
$L_1$ cannot be order embedded into $L_2$. 
The following two sets cannot be separated by a Borel set:

$D_1=\{(L_1,L_2)\in LO^2\;:\;$ $L_1\in WO$ and $L_2\nembed L_1 \}$ 

$D_2=\{(L_1,L_2)\in LO^2\;:\;$ $L_2\in WO$ and $L_1\nembed L_2 \}$ 

\noindent 
To see that these sets are not separable by a Borel set,  first note that for
any $\bPi^1_1$ set $A\sq 2^\om$ there exists a continuous map $f:2^\om\to LO$
such that $f^{-1}(WO)=A$.  (Such a map can be obtained by using the
Kleene-Brouwer ordering on a possible well-founded tree $T\sq \om^{<\om}$ and
mapping $\om^{<\om}\ss T$ to and $\om$ sequence at the end.)  Similar, for any
$\bPi^1_1$ set $B\sq 2^\om$  there exists a continuous map $g:2^\om\to LO$ such
that $g^{-1}(WO)=B$.  Now if $A$ and $B$ happen to be an inseparable disjoint
pair, then the map $h(x)=(f(x),g(x))$ has the property that $h(A)\sq D_1$ and
$h(B)\sq D_2$.  Hence if $C$ separated $D_1$ and $D_2$, then 
$h^{-1}(C)$ would separate $A$ and $B$.

\medskip\noindent Maitra \cite{maitra} uses an open game $G(x)$
on $\om^\om$ due to Blackwell and shows that 

$I= \{x\sq\om^{<\om}\;:\; G(x)$ is won by player I $\}$

$II=\{x\sq\om^{<\om}\;:\; G(x)$ is won by player II $\}$

\noindent  are disjoint inseparable $\bPi^1_1$ sets.  They are not 
complementary sets because in the game considered there may be `ties'.

\medskip\noindent
Becker \cite{becker},\cite{becker2}
contains several examples of inseparable $\bPi^1_1$ sets,
for example,

$B_1=\{f\in C([0,1])\;:\; f \mbox{ is nowhere differentiable }\}$

$B_2=\{f\in C([0,1])\;:\; \exists! x\; f^\prime(x)  \mbox{ exists }\}$

\noindent are inseparable $\bPi^1_1$ sets. He
gives other examples in the compact subsets of the plane:

$C_1=\{K\in {\cal K}({\mathbb R}^2)\::\; K \mbox{ is path-connected and 
simply connected}\}$

$C_2=\{K\in {\cal K}({\mathbb R}^2)\::\; K \mbox{ is path-connected and 
has exactly one hole}\}$

\medskip\noindent 
Milewski \cite{milewski} shows that the following pair of $\bPi^1_1$
sets in the space of compact subsets of the Hilbert
cube, $[0,1]^\om$, are inseparable:

$M_1=\{K\in {\cal K}([0,1]^\om)\;:\;$ all components of $K$ are
finite dimensional $\}$

$M_2=\{K\in {\cal K}([0,1]^\om)\;:\;$ exactly one component of $K$ is
$\infty$-dim $\}$

\medskip\noindent
Camerlo and Darji \cite{cd} give  several families of
pairwise inseparable coanalytic sets.  
For any compact set $K\sq \om^\om$ let

$CD(K)=\{T\sq \om^{<\om}\;:\; \{x\in\om^\om:\forall n \;x\res n\in T\}$
is homeomorphic to $K\}$

\noindent Then for any two nonhomeomorphic compact set $K_1$ and $K_2$ the
sets $CD(K_1)$ and $CD(K_2)$ are inseparable $\bPi^1_1$ sets.

\medskip\noindent
One schema for obtaining natural disjoint inseparable pairs is to 
take a naturally defined filter $F$ on $\om$ and its dual ideal 
$F^*=\{\om\ss X: X\in F\}$.
Note that $F$ and $F^*$ have the same complexity since there exists a 
recursive homeomorphism taking one to other, i.e., $X\mapsto \om\ss X$.
The cofinite filter $COF$ and its dual ideal $FIN$ are naturally inseparable 
$\bSig^0_2$ sets in $P(\om)$.  
Louveau's filter $\ll$ \cite{louveau} 
is an example of a $\bPi^1_1$ filter which
cannot be separated from its dual ideal by a Borel set.  
This filter is on the subsets of $\om^{<\om}$ and is defined as
follows:
$$A\in \ll \mbox{ iff Player I has a winning strategy in the game $J(A)$.}$$
where $J(A)$ is the game:
$$
\begin{array}{llllllll}
\mbox{Player I: } & n_0 &             & n_1  &             & 
n_2 &             & \cdots \\
\mbox{Player II: } &    & m_0\geq n_0 &      & m_1\geq n_1 & 
    & m_2\geq n_2 &  \cdots \\
\end{array}
$$
Player I wins iff for some $k$ all $s\supseteq (m_i:i<k)$
are not in $A$.   This can also be described as follows:
$A\in \ll$ iff $\exists \sigma:\om^{<\om}\to \om\;\forall x\in\om^\om$
if $\forall n\; x(n)\geq \sigma(x\res n)$ then $\exists n\; \forall
s\supseteq x\res n\;\; s\notin A$. 
Although superficially it seems as if
$\ll$ is $\Sigma^1_2$, Louveau proves it is $\Pi^1_1$ by using the
fact that open games are determined and noting that Player I has a 
winning strategy iff Player II does not.  

Louveau proves that any Borel real valued function on a compact
metric space is the $\ll$-limit of a sequence of continuous functions.
Hence $\ll$ is a kind of ultimate generalization of the cofinite filter.
 
\begin{prop}
$\ll$ cannot be separated from its dual
ideal $\ll^*$ by a Borel set.
\end{prop}
\proof
This follows easily from the following in Louveau \cite{louveau}.

\bigskip

Corollaire 8. - Soit $X$ un espace m\'etrisable s\'eparable, $C_1$
et $C_2$ deux parties coanalytiques de $X$.
(ii) Si $C_1$ et $C_2$ sont disjoints, il existe une suite
$(H_u)_{u\in \om^{<\om}}$ de ferm\'es de $X$ telle que 
$$C_1\;\sq\; \liminf_{\ll}H_u\;\sq\; \limsup_{\ll}H_u\;\sq\; X\ss C_2.$$

\bigskip

Recall that  
$$x\in\liminf_{\ll}H_u \;\;\rmiff\;\; \{u \;:\;x\in H_u\}\in \ll$$
and
$$x\in\limsup_{\ll}H_u \;\;\rmiff\;\;\{u \;:\;x\in H_u\}\notin \ll^*.$$

Now take $X=2^\om$ and let $C_1$ and $C_2$ be any two
disjoint inseparable $\bPi^1_1$ sets and take
$H_u\sq 2^\om$ to be the closed sets as in Louveau's Corollaire 8.
Suppose for contradiction
that $B\sq P(\om^{<\om})$ is a Borel set with
$\ll\sq B$ and $\ll^*\cap B=\emp$.  Define
$$Q=\{x\in 2^\om: \{u:x\in H_u\}\in B\}.$$ 
Since $B$ is Borel the set $Q$ is Borel.  Note that
$$\liminf_{\ll}H_u \sq Q \sq \limsup_{\ll}H_u$$
and so $C_1\sq Q$ and $Q\sq 2^\om\ss C_2$ which contradicts
that $C_1$ and $C_2$ cannot be separated.
\qed

\noindent There are plenty of examples of proper $\bPi^1_1$ filters.

$W_1=\{A\sq \om^{<\om}: \neg\exists f\in \om^\om\;\exists^\infty n \;
f\res n\in A\}$

$W_2=\{A\sq \om^{<\om}: \neg\exists f\in \om^\om\;\exists^\infty n \;
\exists s\supseteq f\res n\;\; s\in A\}$

\noindent $W_1$ is the ideal of well-founded subrelations,
$W_2$ is the ideal generated by well-founded subtrees. However,
note that $W_1\sq W_2\sq NWD$ where
$NWD$ is the Borel ideal of nowhere dense subsets of $\om^{<\om}$
defined by

$A\in NWD$ iff $\forall s\; \exists t\supseteq s\;
\forall r\supseteq t\;\; r\notin A$.

\noindent Similarly, 

$W_3=\{A\sq \rationals: A $ is well-ordered $\}$

$W_4=\{A\sq \rationals: cl(A)\sq\rationals $ is compact $\}$

\noindent we have that $W_3\sq W_4\sq NWDQ$ where $NWDQ$ is the Borel
ideal of nowhere dense subsets of the rationals $\rationals$.

Hence, it is the case that each of $W_1,W_2,W_3,W_4$ can be separated
from their duals by a Borel set. 

\bigskip
In Solecki \cite{solecki} it is shown
that for any $\bPi^0_3$ filter $F$ there exists a $\bSig^0_2$ set
$B$ with $F\sq B$ and $F^*\cap B=\emp$.  He leaves open whether
the analogous result holds for $\bPi^0_4$ filters. 
Let $F$ be the cofinite $\times$ cofinite filter on $\om\times\om$, i.e.,
for each $A\sq \om\times\om$ we have that 
$$A\in F \;\;\rmiff\;\; \forall^\infty n \; \forall^\infty  m \;(n,m)\in A\}$$
Then $F$ is a proper $\bSig^0_4$ set (see Kechris 
\cite{kechris} \S 23) and so is its dual ideal 
$F^*$.  In 
Solecki \cite{solecki} Example 1.7,
it is shown that $F$ cannot be separated from $F^*$ by a $\bSig^0_2$ set.
Also according to \cite{solecki} Corollary 1.5, they cannot
be separated by a $\bDel^0_3$ sets.  They can however be separated by
a $\bSig^0_3$ set. Let
$$Q=\{A\sq \om\times\om: \forall^\infty n \; \exists^\infty  m \;(n,m)\in A\}$$
Then $Q$ is $\Sigma^0_3$ and $F\sq Q$ and $F^*\cap Q=\emp$. 

\begin{question}
Is there a $\bSig^0_3$ filter $F$ which cannot be separated from its
dual ideal $F^*$ by a $\bDel^0_3$ set?  In fact, is there a 
$\bSig^0_3$ filter $F$ which is not $\bSig^0_2$?
\end{question}

\begin{question}
For $F$ the cofinite $\times$ cofinite filter does there exist  
a natural $\bSig^0_4$ set $G$ such that $F$ and $G$ are a disjoint
inseparable pair.  (How would you prove there isn't a natural one?)
\end{question}

\bigskip

There is an easy way to generate examples of inseparable
$\bSig^0_n$ sets.  

\begin{prop}
Suppose that $Q\sq 2^\om$ is a complete
$\bPi_n^0$ set.  
Let

$Q_0=\{(x_n:n<\om):\exists n $ even $x_n\in Q$ and 
$\forall m<n\;\; x_m\notin Q\}$ 

$Q_1=\{(x_n:n<\om):\exists n $ odd $x_n\in Q$ and 
$\forall m<n\;\; x_m\notin Q\}$ 

\noindent Then $Q_0$ and $Q_1$ are $\bSig^0_{n+1}$ sets which cannot be
separated by a $\bDel^0_{n+1}$ set. 
\end{prop}

\proof 

Let $A,B\sq 2^\om$ be a disjoint inseparable pair of $\bSig^0_{n+1}$ sets.
Write them as unions of $\bPi^0_n$ sets, 
$A=\cup_{n<\om}U_n^0$ and $B=\cup_{n<\om}U_n^1$.
Since
$Q$ is  complete, there are continuous maps
$f_{2n+i}:2^\om\to 2^\om$ with $f_{2n+i}^{-1}(Q)=U_n^i$.  Then the map
$x\mapsto (f_m(x):m<\om)$ shows that $Q_0$ and $Q_1$ are 
inseparable.
\qed

Similarly there is a natural pair of inseparable $\bSig^0_3$ sets:

\begin{prop}
Let 

$E=\{x\in \om^\om\;:\;\liminf_n\; x(n)$ is even $\}$

$O=\{x\in \om^\om\;:\;\liminf_n\; x(n)$ is odd $\}$

\noindent Then $E$ and $O$ are disjoint inseparable $\bSig^0_3$ sets.
\end{prop}

\proof The set $A=\{x \in \om^\om\;:\;\liminf_n\; x(n)<\infty\}$
is known to be a complete $\bSig^0_3$, see Kechris \cite{kechris}
p.180. This means the given any $\bSig^0_3$ set $B\sq 2^\om$ there
exists a continuous map $f:2^\om\to \om^\om$ with $f(A)=B$.  Now
suppose that $B_1$ and $B_2$ are a disjoint inseparable pair of
$\bSig^0_3$ sets and $f_i$ continuous with $f_i^{-1}(A)=B_i$.  Define
$h:2^\om\to\om^\om$ by $h(x)(n)=2f_1(x(n)$ if $f_1(x)(n)\leq f_2(x)(n)$ 
and
$h(x)(n)= 2 f_2(x(n) +1$ otherwise.  Then $h$ is continuous
and $h(B_1)\sq E$ and $h(B_2)\sq O$ and so $E$ and $O$ cannot
be separated.
\qed

\begin{flushleft}
Arnold W. Miller \\
miller@math.wisc.edu \\
http://www.math.wisc.edu/$\sim$miller\\
University of Wisconsin-Madison \\
Department of Mathematics, Van Vleck Hall \\
480 Lincoln Drive \\
Madison, Wisconsin 53706-1388 \\
\end{flushleft}

\appendix

\newpage
\begin{center}
Appendix A \\ 
\end{center}

This is not intended for publication but only for the electronic version.

\bigskip\noindent Details of the proof of Lemma \ref{mainlemma}.

\bigskip\noindent
{\bf Claim}.
Every nonempty $\Pi^1_1(x)$-set contains a $\Pi^1_1(x)$ singleton.

\proof
Most proofs of the Addison-Kondo Theorem that
every $\Pi^1_1$ set contains a $\Pi^1_1$ singleton
relativizes, e.g., Kechris \cite{kechris}.
It is also follows 
from $\Pi^1_1$ Uniformization property 
(Addison-Kondo Theorem.)
Namely let 
$U\sq \om\times 2^\om\times 2^\om$ be $\Pi^1_1$ set such that for
every $x\in 2^\om$ and
for every set $C$ which is $\Pi^1_1(x)$ there exists $n<\om$ such
that $C=U_(n,x)$.  By the Addison-Kondo Theorem there exists
$V\sq U$ such for every $(n,x)$ if there exists $y$ with
$(n,x,y)\in  U$, then 
there exists a unique $y$ with $(n,x,y)\in V$.)
\qed

\bigskip\noindent
{\bf Claim}. If $(X,R)\in \adm$ is a well-founded relation in
an admissible set $\adm$, then its rank and its canonical rank function 
are in $A$.   
\proof
Define
$\psi(r,D,\al)$ iff 
\begin{enumerate}
\item $r:D\to\al$ is onto the ordinal $\al$, 
\item $D\sq X$, 
\item $\forall x\in D\;\forall  y\in X\;\; 
(yRx\to y\in D)$, and 
\item $\forall x\in D \;\; r(x)=\sup\{r(y)+1:yRx\}$.
\end{enumerate}
Then $\psi$ is a $\Delta_0$ formula.  Also for any $D\sq X$ which
is closed under $R$ both $r$ and $\al$ are unique and this uniqueness
is provable in KP.  Let
$$Q=\{ (r,D,\al)\;:\;
\adm\models \psi(r,D,\al)\}$$
First note that for any  $(r_1,D_1,\al_1), (r_2,D_2,\al_2)\in Q$ that
$$(r_1\cup r_2,D_1\cup D_2,\sup(\al_1,\al_2))\in Q,$$ 
since 
canonical rank functions must agree on their common domain.
Now define $F(x,\be)$ iff there exists $(r,D,\al)\in Q$
with $x\in D$ and $r(x)=\be$.  Then $F$ is $\Sigma_1$ predicate
on $\adm$ which is the graph of
a (possibly partial) function which we also denote $F$.  
By the $\Sigma_1$-replacement axiom of KP there exist
$\de_0\in \adm$ such that for all $x\in X$ and $\beta\in \adm$
$\;\;F(x,\be)\to \be<\de_0$.  First we show that
the domain of $F$ is $X$.  We are assuming that
$(X,R)$ is well-founded, so there exists an $R$-least $x\in X$
such that $x$ is not in the domain of $F$.  Let $R(x)$ be the 
smallest subset of $X$ which contains $\{y\;:\; yRx\}$ and
is closed downward with respect to $R$.  Then
$R(x)\in\adm$ (of course this is obvious
if we assume that $R$ is a strict partial order).   Now
$F\res R(x)\in \adm$ since its graph is a $\Delta_1$ subset of
$R(x)\times \de_0$.  This yields a contradiction since we can then
assign map $x$ to the $\sup\{F(y)+1\;:\; yRx\}$ and get an element
of $Q$ with $x\in D$.  It follows that the domain of $F$ is all
of $X$ and by a similar argument that $F\in\adm$.

\bigskip
\bigskip

Here is a direct proof of the following result of Solecki.

\bigskip\noindent
{\bf Claim}.
Let $F$ be the cofinite $\times$ cofinite filter.
Then $F$ and $F^*$ cannot be separated by
a $\bDel^0_3$ set.

\proof 

First we prove:

\bigskip\noindent
{\bf Lemma}.
Suppose $A$ and $B$ are disjoint $\Sigma^0_3$ subsets of $2^\om$.  Then
there exists a continuous map $h:2^\om\to P(\om\times\om)$ such that
$h(A)\sq F$ and $h(B)\sq F^*$.

\proof
The set 
$$C=\{x\in 2^{\om\times\om} \;:\; \forall^\infty n \; 
\exists  m \;x(n,m)=1\}$$
is a complete $\bSig^0_3$ set, see Kechris 
\cite{kechris} \S 23.  Hence using the theory of Wadge games
there exists a super Lipschitz continuous map $f:2^\om \to 2^{\om\times\om}$
such that $f^{-1}(C)=A$.  By super Lipschitz continuity of $f$ we mean
that $f(x)\res (n\times n)$ is determined by $x\res n$. Let's use
$f^*$ to denote this, i.e., $f(x)\res (n\times n)=f^*(x\res n)$.
The same is true for the set $B$ and let $g$ and $g^*$ be the corresponding
maps. 

Now we use $f^*$ and $g^*$ to construct the map $h^*$ which we think of
as a strategy in a Wadge game.  Fix 
$n_0$.  Given any $s\in 2^{n_0}$ 
assume we have already determined $h^*(s\res (n_0-1))\sq (n_0-1)\times (n_0-1)$.
First of all 
$$h^*(s)\cap (n_0\times (n_0-1)) = h^*(s\res (n_0-1))$$

Given any $n<n_0$ let $i_1^n\leq n$ be the minimal $i$ such that
for all $k$ with $i<k<n$ there exists $m<n_0$ such that
$f^*(s)(n,m)=1$. (If there isn't any such $k$ then
$i_1^n=n$.   Analogously but using $g^*$ define $i_2^n$.  Now
put $(n,n_0)\in h^*(s)$ iff $i_1^n\leq i_2^n$.   In other words, what we are
doing is looking at the n$^{th}$ column and seeing when we look back
at whether $f(x)$ or $g(x)$ is more likely to be in $C$.

The continuous function $h$ is just given by 
$$h(x)=\cup_{n<\om} h^*(x\res n).$$
Now we verify that
$h(A)\sq F$ and $h(B)\sq F^*$. 
Suppose $x\in A$. Since $A$ and $B$ are disjoint we know that
$f(x)\in C$ and $g(x)\notin C$.  This means there exists a $N_0$ so
that for all $n>N_0$ we have that there exists $m$ with $f(x)(n,m)=1$ and
there is $N_1>N_0$ so that $g(x)(N_1,m)=0$ for all $m$.  
(There are infinitely many such columns $N_1$ so just choose the
smallest one bigger than $N_0$.)

We claim
that for all $n>N_1$ the set $h(x)\cap \{n\}\times \om$ is cofinite
in $\{n\}\times \om$.  This is because for a sufficiently large stage
$n_0>n$
in the game the witnesses $m$ will have shown up, i.e. be
less than $n_0$ and so 
$i_1^n$ will be less than or equal to $N_0$ but $i_2^n$
will never be less than $N_1$ and so we will always put $(n,n_0)$
into $h^*(x\res n_0)$.  

The proof that $h(B)\sq F^*$ is analogous.
\qed

The Lemma implies that $F$ and $F^*$ cannot be separated by a $\bDel^0_3$
set, since separation fails for $\bSig^0_3$.

\qed

\newpage
\begin{center}
Appendix B \\ 
\end{center}

\begin{center}
Lecture notes from \\
Slippery Rock conference  \\
Summer Symposium XXVIII June 2004 
\end{center}

For $X$ a Polish space, i.e., separable completely metrizable, define
the Borel classes $\bSig^0_\al$, $\bPi^0_\al$, and 
$\bDel^0_\al$ inductively for countable ordinals $\al$ as follows:

\begin{itemize}
\item $\bSig^0_1$ is the family of open sets in $X$
\item $\bSig^0_\al$ is the family of all countable unions of
sets from $\bigcup_{\be<\al}\bPi_\be^0$
\item $\bPi_\al=\{X\ss A:A\in\bSig^0_\al\}$
\item $\bDel_\al=\bSig^0_\al\cap \bPi^0_\al$
\end{itemize}

The Borel subsets of $X$ are those in $\bigcup_{\al<\om_1}\bSig^0_\al$.
Lebesgue proved that for any uncountable Polish $X$ 
that $\bSig^0_\al\neq \bPi^0_\al$ for any $\al<\om_1$.
For $\Gamma=\bSig^0_\al$ or $\Gamma=\bPi^0_\al$ define the classical
separation principle:

\bigskip
 $\sep(\Gamma)\equiv \forall A,B\in \Gamma\;\; A\cap B=\emp\to
\exists C\in \Delta\;\; A\sq C \rmand A\cap B=\emp$.

\bigskip
\noindent Luzin \cite{luzin} proved that $\sep(\bPi^0_\al)$ holds for 
$1<\al<\om_1$ 
(also $\sep(\bPi^0_1)$ if $X$ is zero dimensional).
He also proved that $\neg\sep(\bSig^0_\al)$.  He gets an inseparable pair
by applying the 
reduction principle to a pair of 
a doubly universal sets, see Kechris \cite{kechris} \S 22 p.171.

The following result answers a question of Dashiell.  It came up
when he was studying Banach spaces of Baire classes of functions
\cite{dashiell} although the question does not appear there.

\setcounter{theorem}{0}

\begin{theorem}\label{o}
Suppose $X$ is a Polish space and $A\sq X$ is
${\bSig}^0_\al$ but not ${\bPi}_\al^0$.  Then
there exists $A^*\sq X$ which is ${\bSig}^0_\al$ such that
$A\cap A^*=\emp$ but there does not exist a ${\bDel}_\al^0$ set $C$  
which separates $A$ and $A^*$, i.e., 
$A\sq C$ and $C\cap A^*=\emp$.
\end{theorem} 

Define $\bSig^1_1$ (or analytic) subsets of $X$ to be the smallest family 
of subsets of $X$ which contains the Borel sets and
is closed under continuous images.  $\bPi^1_1$ is
the family of coanalytic sets or complements of analytic.
Suslin showed that disjoint analytic subsets of $X$ can be separated
by Borel sets and so $\sep(\bSig^1_1)$ holds and
$\bDel_1^1=$Borel.    Luzin's argument
goes thru to show that $\neg\sep(\bPi^1_1)$. 
Dashiell also raised the same question for the coanalytic sets $\bPi^1_1$.
In this case the answer is independent.

\begin{theorem}\label{t}
Suppose $\bPi^1_1$-determinacy holds, then for any $\bPi^1_1$ set $A$
in a Polish space $X$, if
$A$ is not $\bSig^1_1$, then there exists $A^*\sq X$ a $\bPi^1_1$ set 
disjoint from $A$
which cannot be separated from $A$ by a Borel set.
\end{theorem}

\begin{theorem}\label{r}
Suppose $V=L$,  then there exists a $A\sq 2^\om$ and $\bPi^1_1$ set
which is not $\bSig^1_1$ with the
property that for any $B\sq 2^\om$ a $\bPi^1_1$ set disjoint from $A$
there exists a Borel set $C$ with $A\sq C$ and $C\cap B=\emp$.
\end{theorem}

Theorem~\ref{t} is an easy corollary of a result of Steel \cite{steel}:

\begin{theorem}
(Harrington \cite{harrington}, Steel \cite{steel}) The following are
equivalent:

(a) $\bPi^1_1$-determinacy

(b) For all $A\sq X$ and $B\sq Y$ if $X$ and $Y$ are Polish and
$A$ and $B$ are properly analytic, then there exists a Borel
bijection $f:X\to Y$ such that $f(A)=B$.

\end{theorem}
 
Theorem~\ref{r} uses the self-constructible reals $A$ studied
by Guaspari, Sacks, and Kechris, see \cite{kechris2}.
$$A=\{x\in 2^\om\;:\; x\in L_{\om_1^x}\}$$
where $\om_1^x$ is the smallest ordinal not recursive in $x$.

\begin{question}    
If every non Borel $\bPi^1_1$ set is half of an inseparable pair,
then is $\bPi_1^1$-determinacy true?   
\end{question}

Cliff Weil raised the question after the talk of whether we
can get a large number of examples in Theorem \ref{r}, e.g., 

\begin{question}
Assuming V=L,
does there exist continuum many coanalytic sets which are pairwise non Borel
isomorphic and each of which is not half of an inseparable pair?
\end{question}

A number of authors have given natural examples of inseparable pairs
of $\bPi^1_1$ sets, 
Luzin \cite{luzinbook},
Novikov \cite{novikov},
Sierpinski \cite{sier}, 
Dellacherie and Meyer \cite{del},
Maitra \cite{maitra},
Becker \cite{becker},\cite{becker2},
Milewski \cite{milewski},
and Camerlo and Darji \cite{cd}.

Another method
for obtaining a disjoint inseparable pair is to 
take a filter $F$ on $\om$ and its dual ideal  $F^*=\{\om\ss X: X\in F\}$. 
Note that $F$ and $F^*$ have the same complexity since there exists a 
recursive homeomorphism taking one to other, i.e., $X\mapsto \om\ss X$.
This was suggested by the results in Solecki \cite{solecki}.

The cofinite filter $COF$ and its dual ideal $FIN$ are naturally inseparable 
$\bSig^0_2$ sets in $P(\om)$.  
Louveau's filter $\ll$ \cite{louveau} 
is an example of a proper $\bPi^1_1$ filter. 
Louveau proves that the Borel real valued function on a compact
metric space are exactly the $\ll$-limits of sequences of continuous 
functions.
Hence $\ll$ is a kind of ultimate generalization of the cofinite filter.
 
\begin{prop}
$\ll$ cannot be separated from its dual
ideal $\ll^*$ by a Borel set.
\end{prop}

\end{document}